\author{\Large{Damanvir Singh Binner}}
\documentclass[12pt]{article}
\usepackage[margin=1in]{geometry}
\usepackage[usenames]{xcolor}
\usepackage{amssymb}
\usepackage{relsize}
\usepackage{amsmath}
\usepackage{amsthm}
\usepackage{amsfonts}
\usepackage{amscd}
\usepackage{graphicx}
\usepackage{hyperref}
\usepackage{tikz-cd}
\usepackage{tikz}

\begin{document}

\newcommand{\D}[1]{{\bf \color{red} #1}}
\newcommand{\M}[1]{{\bf \color{magenta} #1}}
\newcommand{\K}[1]{{\bf \color{violet} #1}}

\theoremstyle{plain}
\newtheorem{theorem}{Theorem}
\newtheorem{corollary}[theorem]{Corollary}
\newtheorem{lemma}[theorem]{Lemma}
\newtheorem{proposition}[theorem]{Proposition}
\newtheorem{question}[theorem]{Question}

\theoremstyle{definition}
\newtheorem{definition}[theorem]{Definition}
\newtheorem{example}[theorem]{Example}
\newtheorem{conjecture}[theorem]{Conjecture}

\theoremstyle{remark}
\newtheorem{remark}[theorem]{Remark}

\title{\Large{Number of partitions of $n$ with a given parity of the smallest part}}
\date{}
\maketitle
\begin{center}
\vspace*{-8mm}
\large{Department of Mathematics \\
Indian Institute of Science Education and Research (IISER) \\
Mohali, Punjab, India \\
damanvirbinnar@iisermohali.ac.in}
\end{center}

\begin{abstract}
We obtain a combinatorial proof of a surprising weighted partition equality of Berkovich and Uncu. Our proof naturally leads to a formula for the number of partitions of $n$ with a given parity of the smallest part, in terms of $S(i)$, the number of partitions of $i$ into distinct parts with even rank minus the number with odd rank, for which there is an almost closed formula by Andrews, Dyson and Hickerson. This method of calculating the number of partitions of $n$ with a given parity of the smallest part is practical and efficient.
\end{abstract}

\section{Introduction}
\label{S1}

In recent times, the properties of partitions based on the parity of their parts have received special attention (see \cite{And1, And2, And3, Saikia, Lovejoy}). In this direction, Berkovich and Uncu \cite{berkbook} proved a rather surprising result connecting two very simple looking partition statistics. For a given partition $\pi$ of $n$ (denoted by $\pi \vdash n$), let $s(\pi)$ denote the smallest part of $\pi$, let $f_i$ denote the number of times $i$ appears as a part in $\pi$, and let $t(\pi)$ denote the number defined by the following properties. 

\begin{enumerate}
\item $f_i \equiv 1$ (mod $2$) for $1 \leq i \leq t(\pi)$.
\item $f_{t(\pi)+1} \equiv 0$ (mod $2$).
\end{enumerate}

The authors referred to $t(\pi)$ as the \emph{length of the initial odd-frequency chain}.   

\begin{theorem}[Berkovich and Uncu (2017)]
\label{Cute}
For any $n \geq 1$, $$ \sum_{\pi \vdash n} (-1)^{s(\pi)+1} = \sum_{\pi \vdash n} t(\pi). $$ 
\end{theorem}

Berkovich and Uncu proved Theorem \ref{Cute} in \cite[Theorem 3.1]{berkbook} using Jackson's transformation, and also gave an alternate $q$-series proof in \cite[Section 5]{BU19}. See \cite[Table 1]{BU19} for an illustration of Theorem \ref{Cute}. However, this elegant result definitely calls for an elementary combinatorial proof. This is the main goal of this note. We prove Theorem \ref{Cute} by combinatorially proving Theorems \ref{Smallest} and \ref{OF}, described below. Prior to that, we also need to introduce the quantity $S(i)$, the number of partitions of $i$ into distinct parts with even rank minus the number with odd rank. Andrews, Dyson and Hickerson \cite[Theorem 2, Theorem 3]{AD} gave an almost closed formula for $S(i)$ described in Theorems \ref{StoT} and \ref{Powerful} below. 

\begin{theorem}[Andrews, Dyson and Hickerson (1988)]
\label{StoT}
For $n \geq 0$, $S(n) = T(24n+1)$, where $T(m)$ is as defined in \cite[Section 1]{AD}.
\end{theorem}

\begin{theorem}[Andrews, Dyson and Hickerson (1988)]
\label{Powerful}
Let $m \neq 1$ be an integer $\equiv 1$ (mod $6$). Suppose $m$ has the factorization 
\begin{equation}
\label{Factor}
m=p_1^{e_1}p_2^{e_2} \cdots p_r^{e_r},
\end{equation}
 where $r \geq 1$, each $p_i$ is either a prime $\equiv 1$ (mod $6$) or negative of a prime $\equiv 5$ (mod $6$), the $p_i$'s are distinct, and the $e_i$'s are positive integers. Then, $T(m) = T(p_1^{e_1}) \cdots T(p_r^{e_r})$, where 
$$ 
T(p^e) = 
\begin{cases}
0, & \text{if } p \not\equiv 1 (\bmod 24) \text{and $e$ is odd}, \\
1, & \text{if } \text{$p \equiv 13$ or $19$} (\bmod 24) \text{and $e$ is even}, \\
(-1)^{\frac{e}{2}}, & \text{if } p \equiv 7 (\bmod 24) \text{and $e$ is even}, \\
e+1, & \text{if } p \equiv 1 (\bmod 24) \text{and $T(p)=2$}, \\
(-1)^e (e+1), & \text{if } p \equiv 1 (\bmod 24) \text{and $T(p)=-2$},
\end{cases}
$$

\end{theorem}
We say that Theorem \ref{Powerful} gives an almost closed formula for finding $T(n)$, because it does not determine the sign of $T(p)$ when $p \equiv 1$ (mod $24$). We will describe how to deal with this issue later.  Let $P_O(n)$ (respectively $P_E(n)$) denote the number of partitions of $n$ whose smallest part is odd (respectively even). Then, Theorem \ref{Cute} can be rewritten as 
\begin{equation}
\label{Alt}
P_O(n)-P_E(n) = \sum_{\pi \vdash n} t(\pi). 
\end{equation}

\begin{theorem}
\label{Smallest}
For any $n \geq 1$, $$P_O(n)-P_E(n) = \sum_{i \geq 1} S(i)p(n-i). $$
\end{theorem}

\begin{remark}
Though Berkovich and Uncu did not mention Theorem \ref{Smallest} or its $q$-series version, the latter can be easily obtained from their alternate proof of Theorem \ref{Cute} in \cite[Section 5]{BU19}. As described below, Theorem \ref{Smallest} leads to a practical and efficient method to calculate the number of partitions of $n$ of a given parity, and is naturally obtained using our combinatorial approach. 
\end{remark}

\begin{theorem}
\label{OF}
For any $n \geq 1$, $$ \sum_{\pi \vdash n} t(\pi) =  \sum_{i \geq 1} S(i)p(n-i). $$ 
\end{theorem}

We describe how Theorem \ref{Smallest} can be used to calculate the number of partitions of $n$ with odd and even smallest parts. Note that the calculation of $p(n)$ is very easy compared to finding all the partitions of $n$ because of the availability of formulae such as Hardy-Ramanujan-Rademacher formula and recurrences such as Euler's recurrence. The calculation of $S(i)$ is also relatively very easy because of Theorem \ref{Powerful}. As mentioned before, the biggest difficulty that one faces is finding the sign of $T(p)$ when $p \equiv 1$ (mod $24$). We briefly describe the method to calculate the sign. For details, refer to \cite[pp.\ $398$--$399$]{AD}. Suppose $p \equiv 1$ (mod $24$) appears in the factorization of $24i+1$ in the form of \eqref{Factor}. From the standard theory of generalized Pell's equation \cite[Theorem 3.3]{Conrad2}, there exists a solution of $x^2 - 6y^2=p$ such that $0 \leq y < \sqrt{|p|}$. We find the least nonnegative integer $y_0$ such that $6y_0^2+p$ is a perfect square, say $x_0^2$ for some nonnegative integer $x_0$. Then, $T(p)$ is positive if $x_0+3y_0 \equiv \pm 1$ (mod $12$) and negative if $x_0+3y_0 \equiv \pm 5$ (mod $12$).

We demonstrate Theorem \ref{Smallest} for $n=37$. That is, we find the number of partitions of $37$ that have an odd smallest part, and the number of partitions of $37$ that have an even smallest part. The values of $S(i)$ can be calculated using Theorem \ref{Powerful} and the strategy in the above paragraph. Most of the cases follow directly from Theorem \ref{Powerful}. The values of $S(i)$ for these cases are recorded in Table \ref{tab1}. The cases which need more work are $i=3,4,8,10,13,14,17,18,19,24,25,28,32$. In these cases, $24i+1$ is a prime number. For example, suppose $i=3$. We need to find the least nonnegative integer $y_0$ such that $6y_0^2+73$ is a perfect square. It is easy to see that $y_0 = 4$. Then $x_0 = 13$, and thus $x_0+3y_0 = 25 \equiv 1$ (mod $12$). Therefore, $S(3) = T(73)$ is positive. The values of $x_0$, $y_0$ and $S(i)$ for these cases are recorded in Table \ref{tab1}. 
\begin{table}[htpb]
    \centering
    \begin{tabular}{|c|c|c|c|c|c|}
 \hline
$i$      &   $24i+1$ & $y_0$ & $x_0$ & $x_0 + 3y_0$ (mod $12$) & $S(i) = T(24i+1)$ \\
\hline
$1$ & $25$ & - & - & - & $1$ \\
\hline
$2$ & $49$ & - & - & - & $-1$\\
\hline
$3$ & $73$ & $4$ & $13$ & $1$ & $2$  \\
\hline
$4$ & $97$ & $2$ & $11$ & $5$ & $-2$ \\
\hline
$5$ & $121$ & - & - & - & $1$ \\
\hline
$6$ & $145$ & - & - & - & $0$ \\
\hline
$7$ & $169$ & - & - & - & $1$ \\
\hline
$8$ & $193$ & $4$ & $17$ & $5$ & $-2$ \\
\hline
$9$ & $217$ & - & - & - & $0$ \\
\hline
$10$ & $241$ & $8$ & $25$ & $1$ & $2$ \\
\hline
$11$ & $265$ & - & - & - & $0$ \\
\hline
$12$ & $289$ & - & - & - & $-1$ \\
\hline
$13$ & $313$ & $6$ & $23$ & $5$ & $-2$ \\
\hline
$14$ & $337$ & $2$ & $19$ & $1$ & $2$ \\
\hline
$15$ & $361$ & - & - & - & $1$ \\
\hline
$16$ & $385$ & - & - & - & $0$ \\
\hline
$17$ & $409$ & $6$ & $25$ & $-5$ & $-2$ \\
\hline
$18$ & $433$ & $4$ & $23$ & $-1$ & $2$ \\
\hline
$19$ & $457$ & $8$ & $29$ & $5$ & $-2$\\
\hline
$20$ & $481$ & - & - & - & $0$\\
\hline
$21$ & $505$ & - & - & - & $0$ \\
\hline
$22$ & - & - & - & - & $3$ \\
\hline
$23$ & $553$ & - & - & - & $0$\\
\hline
$24$ & $577$ & $8$ & $31$ & $-5$ & $-2$ \\
\hline
$25$ & $601$ & $2$ & $25$ & $-5$ & $-2$ \\
\hline
$26$ & $625$ & - & - & - & $1$ \\
\hline
$27$ & $649$ & - & - & - & $0$ \\
\hline
$28$ & $673$ & $14$ & $43$ & $1$ & $2$ \\
\hline
$29$ & $697$ & - & - & - & $0$ \\
\hline
$30$ & $721$ & - & - & - & $0$ \\
\hline
$31$ & $745$ & - & - & - & $0$ \\
\hline
$32$ & $769$ & $10$ & $37$ & $-5$ & $-2$ \\
\hline
$33$ & $793$ & - & - & - & $0$ \\
\hline
$34$ & $817$ & - & - & - & $0$ \\
\hline
$35$ & $841$ & - & - & - & $1$ \\
\hline
$36$ & $865$ & - & - & - & $0$ \\
\hline
$37$ & $889$ & - & - & - & $0$ \\
\hline
    \end{tabular}
    \caption{The values of $S(i)$ for $1 \leq i \leq 37$}
        \label{tab1}
\end{table}

Thus, using Theorem \ref{Smallest} and the vaules of $S(i)$ given in Table \ref{tab1}, we have 
\begin{align*}
 P_O(37) - P_E(37) &= p(36)-p(35)+2p(34)-2p(33)+p(32)+p(30)-2p(29)+2p(27) \\
 &-p(25)-2p(24)+2p(23)+p(22)-2p(20)+2p(19)-2p(18) \\
 &+3p(15) -2p(13)-2p(12)+p(11)+2p(9)-2p(5)+p(2). 
 \end{align*}
 
 Substituting the values, we get that $$P_O(37)-P_E(37) = 15907.$$ Moreover, we have $$P_O(37) + P_E(37) = p(37) = 21637.$$ From these two equations, it follows that $ P_O(37) = 18772$ and $ P_E(37) = 2865$. That is, there are $18772$ partitions of $37$ that have an odd smallest part, and $2865$ partitions of $37$ that have an even smallest part. 
 
 In fact, using the vaules of $S(i)$ given in Table \ref{tab1}, we can find the number of partitions of $n$ with a given parity for any $n \leq 37$. For example, suppose $n=17$.  Using Theorem \ref{Smallest} and the vaules of $S(i)$ given in Table \ref{tab1}, we have 
\begin{align*}
 P_O(17) - P_E(17) &= p(16)-p(15)+2p(14)-2p(13)+p(12)+p(10)-2p(9)+2p(7) \\
 &-p(5)-2p(4)+2p(3)+p(2)-2p(0). 
 \end{align*}
 
 Substituting the values, we get that $$P_O(17)-P_E(17) = 201.$$ Moreover, we have $$P_O(17) + P_E(17) = p(17) = 297.$$ From these two equations, it follows that $ P_O(17) = 249 $ and $ P_E(17) = 48 $. That is, there are $249$ partitions of $17$ that have an odd smallest part, and $48$ partitions of $17$ that have an even smallest part.

\section{Proof of Theorem \ref{Smallest}}

We calculate the number of partitions with a given smallest part $i$ using the principle of inclusion and exclusion (PIE). This approach was also used by the present author and Rattan in the author's PhD Thesis \cite[Section 5.1]{Thesis} to obtain a natural combinatorial proof of Euler's recurrence for integer partitions using PIE. Since there are crucial differences, we provide all the details here for the sake of completeness. We recall some notation defined in \cite[Section 5.1]{Thesis} and define some new notation. 

\begin{itemize}
\item $A_{j,k}(n)$ is the set of partitions of $n$ having exactly $k$ parts of size $j$;
\item $B_{j,k}(n)$ is the set of partitions of $n$ having at least $k$ parts of size $j$. 
\end{itemize}

The following properties of these sets are immediate. 

\begin{enumerate}
\item $|B_{j,k}(n)| = p(n-jk)$.
\item If $j \neq j'$, then $|B_{j,k}(n) \cap B_{j',k'}(n)| = p(n-jk-j'k')$.
\item $A_{j,0}(n) = B_{j,1}^{\complement}(n)$, where the complementation is with respect to the set $Par(n)$, consisting of all partitions of $n$. 
\end{enumerate}

We will also need the following sets.

\begin{itemize}
\item $T_{s,i}$ denotes the set of partitions into $s$ distinct parts with largest part less than $i$. 
\item $U_{s,i}$ denotes the set of partitions into $s$ distinct parts with largest part equal to $i$. 
\end{itemize}

Recall that $r(\pi)$ denotes the rank of $\pi$, which is obtained by subtracting the number of parts of $\pi$ from the largest part of $\pi$. Further, let $D$ denote the set of all nonempty partitions into distinct parts. 

Clearly, the number of partitions of $n$ with smallest part $i$ is equal to the number of partitions of $n-i$ with no part equal to $1,2,\ldots,(i-1)$. Then, using PIE along with the definitions and properties of the sets $A_{j,k}$ and $B_{j,k}$, it follows that the number of partitions of $n$ with smallest part $i$ is given by

\begin{align*} 
 &|A_{1,0}(n-i) \cap A_{2,0}(n-i) \cdots \cap A_{i-1,0}(n-i)| \\
 =& |B_{1,1}^{\complement}(n-i) \cap B_{2,1}^{\complement}(n-i) \cdots \cap B_{i-1,1}^{\complement}(n-i)| \\
 =&|(B_{1,1}(n-i) \cup B_{2,1}(n-i) \cdots \cup B_{i-1,1}(n-i))^{\complement}| \\
 =& \sum_{s=0}^{i-1} (-1)^s \sum_{(i_1, i_2, \ldots, i_s) \in T_{s,i}} p(n-i-i_1-i_2 \cdots -i_s) \\
 =& \sum_{s=0}^{i-1} (-1)^s \sum_{\pi \in U_{s+1,i}} p(n-|\pi|) \\
 =& \sum_{s=1}^{i} (-1)^{s-1} \sum_{\pi \in U_{s,i}} p(n-|\pi|). 
\end{align*}
Therefore, 
\begin{align*} 
P_O(n)-P_E(n) = & \sum_{i \geq 1} (-1)^{i+1} \sum_{s=1}^{i} (-1)^{s-1} \sum_{\pi \in U_{s,i}} p(n-|\pi|) \\
&=\sum_{i \geq 1} \sum_{s=1}^{i}  \sum_{\pi \in U_{s,i}} (-1)^{i-s} p(n-|\pi|) \\
&=\sum_{i \geq 1} \sum_{s=1}^{i}  \sum_{\pi \in U_{s,i}} (-1)^{r(\pi)} p(n-|\pi|) \\
&=\sum_{\pi \in D} (-1)^{r(\pi)} p(n-|\pi|) \\
&= \sum_{i \geq 1} S(i)p(n-i),
\end{align*}
as required.

\section{Proof of Theorem \ref{OF}}

We define a partition $\pi$ to be a \emph{$C$-partition} if for all $i \geq 1$, whenever $i$ appears as a part in $\pi$, then all natural numbers less than $i$ must also appear as a part in $\pi$. We let $V$ denote the set of all nonempty $C$-partitions, and let $V_{s,i}$ denote the set of all $C$-partitions with largest part $i$ such that $s$ of the numbers from $1$ to $i$ appear with an even frequency. 

Further, let $W_i$ denote the set of partitions in which the numbers $1,2, \ldots, i$ appear with an odd frequency. From the definition of $t(\pi)$, it follows that $$ \sum_{\pi \vdash n} t(\pi) = \sum_{i \geq 1} |W_i|. $$ We again use PIE to calculate $|W_i|$.
Clearly, the number of partitions of $n$ with $m_j$ parts of $j$ for all $1 \leq j \leq i$ is equal to the number of partitions of $n-\sum_{j=1} jm_j$ that have no parts of $1,2,\ldots,i$. Therefore, using PIE along with the definitions and properties of the sets $A_{j,k}$ and $B_{j,k}$, we have 

\begin{equation}
\label{OFCL}
\begin{aligned} 
&\sum_{\pi \vdash n} t(\pi) \\
&= \sum_{i \geq 1} |W_i|  \\
 &= \sum_{i \geq 1} \sum_{\substack{m_1,\ldots,m_i \\ odd}} \left|A_{1,0}\left(n-\sum_{j=1}^i jm_j\right) \cap A_{2,0}\left(n-\sum_{j=1}^i jm_j\right) \cdots \cap A_{i,0}\left(n-\sum_{j=1}^i jm_j\right)\right|  \\
 &= \sum_{i \geq 1} \sum_{\substack{m_1,\ldots,m_i \\ odd}} \left|B_{1,1}^{\complement}\left(n-\sum_{j=1}^i jm_j\right) \cap B_{2,1}^{\complement}\left(n-\sum_{j=1}^i jm_j\right) \cdots \cap B_{i,1}^{\complement}\left(n-\sum_{j=1}^i jm_j\right)\right|  \\
 &= \sum_{i \geq 1} \sum_{\substack{m_1,\ldots,m_i \\ odd}} \left|\left(B_{1,1}\left(n-\sum_{j=1}^i jm_j\right) \cup B_{2,1}\left(n-\sum_{j=1}^i jm_j\right) \cdots \cup B_{i,1}\left(n-\sum_{j=1}^i jm_j\right)\right)^{\complement}\right|  \\
&= \sum_{i \geq 1} \sum_{\substack{m_1,\ldots,m_i \\ odd}} \sum_{s=0}^{i} (-1)^s \sum_{(i_1, i_2, \ldots, i_s) \in T_{s,i}} p\left(n-\sum_{j=1}^i jm_j -(i_1 + \cdots + i_s) \right)  \\
&= \sum_{i \geq 1}  \sum_{s=0}^{i} \sum_{\substack{m_1,\ldots,m_i \\ odd}} \sum_{(i_1, i_2, \ldots, i_s) \in T_{s,i}} (-1)^s p\left(n-\sum_{j=1}^i jm_j -(i_1 + \cdots + i_s) \right)  \\
&= \sum_{i \geq 1}  \sum_{s=0}^{i} \sum_{\pi \in V_{s,i}} p(n-|\pi|) \\
&= \sum_{\pi \in V} (-1)^{h(\pi)} p(n-|\pi|),
\end{aligned}
\end{equation}

where $h(\pi)$ denotes the number of parts of $\pi$ that have an even frequency.

Using Ferrers diagram, one can show that the conjugation map provides a bijection between the set $V$ of nonempty $C$-partitions and the set $D$ of nonempty partitions into distinct parts such that the $C$-partitions in $V$ with odd $h(\pi)$ (respectively even $h(\pi)$) are mapped to partitions into distinct parts in $D$ with odd rank (respectively even rank). We describe this in some detail. 

We need to simultaneously work with two different notations for a given partition $\pi$. First is the standard notation $\pi = (\pi_1,\pi_2, \dots)$ where $\pi_1 \geq \pi_2 \geq \cdots$. We call this as notation A for $\pi$. In notation B, we write $\pi = (1^{f_1}, 2^{f_2},\ldots)$, where $f_i$ is the \emph{frequency} of $i$ or the number of times a part $i$ occurs in $\pi$. Clearly $f_i \geq 0$ for all $i$. Further, note that for a $C$-partition $\pi$ with largest part $m$, $f_i \geq 1$ for all $1 \leq i \leq m$. 

We know that the conjugation map is an involution on the set $Par(n)$. An algebraic description of the conjugation map can be easily written as follows. Suppose we have a partition $\pi$ which can be expressed in notation B as $$\pi = (1^{f_1}, 2^{f_2},\ldots, k^{f_k})$$ for some $k \in \mathbb{N}$. Then, the image $\bar{\pi}$ of $\pi$ under the conjugation map can be expressed in notation A as 
\begin{equation}
\label{Conju}
 \bar{\pi} = \left(\sum_{i=1}^k f_i, \sum_{i=2}^k f_i, \ldots, f_{k-1}+f_k, f_k \right). 
 \end{equation}

For example, consider the partition $\pi = (7, 7, 6, 4, 4)$ of $28$ (notation A). Clearly, it can be expressed as $(1^0,2^0,3^0,4^2,5^0,6^1,7^2)$ in notation B. The Ferrers diagram for $\pi$ is given as

 
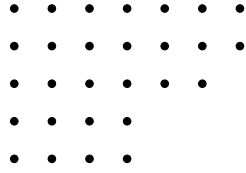
\vspace{.3cm}
 
 \begin{center}
\begin{tikzpicture}

\foreach \x in {0,...,6}
	\filldraw (\x*.5, -5) circle (.5mm);
\foreach \x in {0,...,6}
	\filldraw (\x*.5, -5.5) circle (.5mm);
\foreach \x in {0,...,5}
	\filldraw (\x*.5, -6) circle (.5mm);
\foreach \x in {0,...,3}
	\filldraw (\x*.5, -6.5) circle (.5mm);
\foreach \x in {0,...,3}
	\filldraw (\x*.5, -7) circle (.5mm);.
\end{tikzpicture}
\end{center}

 \vspace{.3cm}

Then, the Ferrers diagram for the conjugate partition $\bar{\pi}$ is given as 

 
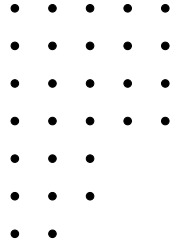
\vspace{.3cm}

\begin{center}
\begin{tikzpicture}

\foreach \x in {0,...,4}
	\filldraw (\x*.5, -5) circle (.5mm);
\foreach \x in {0,...,4}
	\filldraw (\x*.5, -5.5) circle (.5mm);
\foreach \x in {0,...,4}
	\filldraw (\x*.5, -6) circle (.5mm);
\foreach \x in {0,...,4}
	\filldraw (\x*.5, -6.5) circle (.5mm);
\foreach \x in {0,...,2}
	\filldraw (\x*.5, -7) circle (.5mm);
\foreach \x in {0,...,2}
	\filldraw (\x*.5, -7.5) circle (.5mm);
\foreach \x in {0,...,1}
	\filldraw (\x*.5, -8) circle (.5mm);
\end{tikzpicture}
\end{center}

The diagram shows that the partition $\bar{\pi}$ can be expressed as $(5,5,5,5,3,3,2)$ in notation A, verifying the description of the conjugation map provided in \eqref{Conju}.

From this description, it immediately follows that under the conjugation map, the set of all nonempty $C$-partitions of $n$ is mapped to the set of nonempty partitions of $n$ into distinct parts. Finally, we observe that 

\begin{align*} 
r(\bar{\pi}) &= \sum_{i=1}^k f_i - k \\
&= \sum_{i=1}^k \left(f_i - 1\right) \\
&\equiv h(\pi) \pmod 2.
\end{align*}

Thus, using \eqref{OFCL} and the above properties of the conjugation map, we have
\begin{align*} 
\sum_{\pi \vdash n} t(\pi) &= \sum_{\pi \in V} (-1)^{h(\pi)} p(n-|\pi|) \\
&= \sum_{\pi \in D} (-1)^{r(\pi)} p(n-|\pi|) \\
&= \sum_{i \geq 1} S(i)p(n-i),
\end{align*}
as required.

\section{Acknowledgement}
 
 The author acknowledges the support of IISER Mohali for providing research facilities and fellowship.

\section{Data availability}

All data generated or analysed during this study are included in this article.

\end{document}